\documentclass[12pt]{amsart}
\usepackage{amssymb}
\usepackage{epic}
\usepackage{eepic}
\newcommand{\w}{\omega}
\newcommand{\G}{\Gamma}
\newcommand{\ch}{{\C H}}
\newcommand{\wbar}{\bar{\w}}
\newcommand{\thetabar}{\bar{\theta}}
\newcommand{\ip}[2]{\langle #1|#2\rangle}
\newcommand{\spanof}[1]{\langle #1\rangle}
\newcommand{\tensor}{\otimes}
\DeclareMathOperator{\Realpart}{{\rm Re}}
\renewcommand{\Re}{\Realpart}
\DeclareMathOperator{\aut}{{\rm Aut}}
\DeclareMathOperator{\diag}{{\rm diag}}
\newcommand{\code}{\mathcal{C}}
\newcommand{\codeperp}{\code^{\perp}}
\newcommand{\SL}{{\rm SL}}
\newcommand{\cll}{\Lambda_{24}^\E}
\newcommand{\suz}{{\it Suz}}
\newcommand{\Co}{{\it Co}}

\newcommand{\D}{\Delta}
\newcommand{\p}{\rho}
\newcommand{\pperp}{\p^{\perp}}
\renewcommand{\d}{\delta}
\newcommand{\normalin}{\vartriangleleft}
\newcommand{\E}{\mathcal{E}}
\newcommand{\C}{\mathbb{C}}
\newcommand{\F}{\mathbb{F}}
\newcommand{\Z}{\mathbb{Z}}
\newcommand{\Sp}{{\rm Sp}}
\let\iso\cong
\let\cong\equiv
\newcommand{\sset}{\subseteq}

\newtheorem{theorem}{Theorem}
\newtheorem{lemma}[theorem]{Lemma}
\theoremstyle{remark}
\newtheorem*{remark}{Remark}
\newtheorem*{remarks}{Remarks}

\begin{document}
\title{On the $Y_{555}$ complex reflection group}
\author{Daniel Allcock}
\address{Department of Mathematics\\University of Texas, Austin}
\email{allcock@math.utexas.edu}
\urladdr{http://www.math.utexas.edu/\textasciitilde allcock}
\thanks{Partly supported by NSF grant DMS-0600112.}
\date{February 6, 2008}
%\date{January 15, 2008}
%\date{January 4, 2008}

\begin{abstract}
We give a computer-free proof of a theorem of Basak, describing the
group generated by 16 complex reflections of order~3, satisfying the
braid and commutation relations of the $Y_{555}$ diagram.  The group
is the full isometry group of a certain lattice of signature $(13,1)$
over the Eisenstein integers $\Z[\sqrt[3]{1}]$.  Along the way we
enumerate the cusps of  this lattice and
classify the root and Niemeier lattices over $\Z[\sqrt[3]{1}]$
\end{abstract}

\maketitle

The author has conjectured \cite{monstrous} that the largest sporadic finite
simple group, the monster, is related to complex algebraic geometry,
with a certain complex hyperbolic orbifold acting as a sort of
intermediary.  In particular, the bimonster $(M\times M){:}2$ and a
certain group $P\G$ acting on complex hyperbolic 13-space are
conjecturally both quotients of $\pi_1\bigl((\ch^{13}-\D)/P\G\bigr)$
for a certain hyperplane arrangement $\D$ in $\ch^{13}$, got by
adjoining very simple relations.  If the conjecture is true then it
has the consequence that $P\G$ is generated by 16 complex reflections
of order 3, satisfying the braid and commutation relations of the
$Y_{555}$ diagram
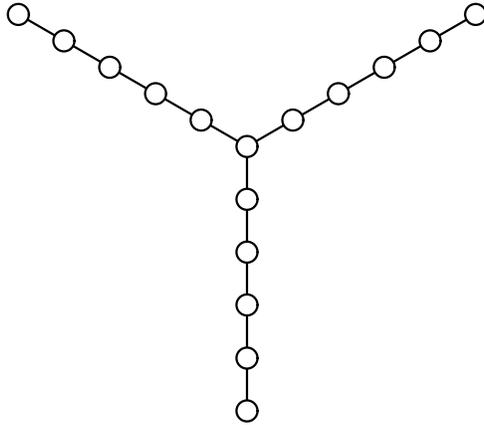
\begin{displaymath}
%
% the suffix D is a reminder that these variables should be dimens
\newdimen\edgelengthD  \edgelengthD=20pt
\newdimen\nodediameterD \nodediameterD=8pt
\newdimen\XshiftD \XshiftD=\edgelengthD \multiply\XshiftD by 173
\divide\XshiftD by 200
\newdimen\YshiftD \YshiftD=\edgelengthD \divide\YshiftD by 2
%
%
% U and V are points to use as local variables
\newcount\UxC
\newcount\UyC
\newcount\VxC
\newcount\VyC
\newcount\nodediameterC
% from (#1,#2) to (#3,#4) arguments should be dimensions
\def\bond#1#2#3#4{%
  \thicklines
  \UxC=#1
  \UyC=#2
  \VxC=#3
  \VyC=#4
  \drawline(\UxC,\UyC)(\VxC,\VyC)
}
% #1,#2 should be an ordered pair of dimensions
\def\hollownode#1#2{%
  \thicklines
  \UxC=#1
  \UyC=#2
  \nodediameterC=\nodediameterD
  \filltype{white}\put(\UxC,\UyC){\circle*{\nodediameterC}}
}
%
% FOR PLACING TEXT AND OTHER STUFF
%
% like ordinary put, except that args should be dimens
\def\myput#1#2#3{%
  \UxC=#1
  \UyC=#2
  \put(\UxC,\UyC){#3}
}
% (#1,#2) should be position (dimensions not numbers)
% (#3,#4) should be how many hundredths of a noderadius
% to offset the position.  (#5,#6) should be an additional
% offset (again, a dimension).  #7 should be the positioning
% string (tr, b, etc.)#8 should the the thing to typset.
%
\newcount\XoffC
\newcount\YoffC
\newdimen\XoffD
\newdimen\YoffD
\newdimen\noderadiusD
\newdimen\UxD
\newdimen\UyD
\def\nearnode#1#2#3#4#5#6#7#8{%
  \UxD=#1
  \UyD=#2
  \XoffC=#3
  \YoffC=#4
  \noderadiusD=\nodediameterD \divide\noderadiusD by 2
  \XoffD=\noderadiusD
  \YoffD=\noderadiusD
  \multiply\XoffD by \XoffC
  \multiply\YoffD by \YoffC
  \divide\XoffD by 100
  \divide\YoffD by 100
  \advance\UxD by \XoffD
  \advance\UyD by \YoffD
  \myput\UxD\UyD{\kern#5\makebox(0,0)[#7]{\raise#6\hbox{#8}}}
}
\setlength{\unitlength}{1sp}
\newcount\LLy
\newcount\height
\LLy=-\edgelengthD
\multiply\LLy by 5
\height=\edgelengthD
\multiply\height by 15
\divide\height by 2
\begin{picture}(0,\height)(0,\LLy)
%\put(0,\LLy){\framebox(0,\height)[bl]{\relax}}
% origin at A
%
%  Fone					Ftwo
%      Eone			    Etwo
%	   Done			Dtwo
%	       Cone	    Ctwo
%		   Bone	Btwo
%		       A
%		     Bthree
%
%		     Cthree
%
%		     Dthree
%
%		     Ethree
%
%		     Fthree
%
\newdimen\AX \AX=0pt
\newdimen\AY \AY=0pt
\newdimen\BoneX \BoneX=-\XshiftD
\newdimen\BoneY \BoneY=\YshiftD
\newdimen\ConeX \ConeX=-\XshiftD \multiply\ConeX by 2
\newdimen\ConeY \ConeY=\YshiftD \multiply\ConeY by 2
\newdimen\DoneX \DoneX=-\XshiftD \multiply\DoneX by 3
\newdimen\DoneY \DoneY=\YshiftD \multiply\DoneY by 3
\newdimen\EoneX \EoneX=-\XshiftD \multiply\EoneX by 4
\newdimen\EoneY \EoneY=\YshiftD \multiply\EoneY by 4
\newdimen\FoneX \FoneX=-\XshiftD \multiply\FoneX by 5
\newdimen\FoneY \FoneY=\YshiftD \multiply\FoneY by 5
\newdimen\BtwoX \BtwoX=\XshiftD
\newdimen\BtwoY \BtwoY=\YshiftD
\newdimen\CtwoX \CtwoX=\XshiftD \multiply\CtwoX by 2
\newdimen\CtwoY \CtwoY=\YshiftD \multiply\CtwoY by 2
\newdimen\DtwoX \DtwoX=\XshiftD \multiply\DtwoX by 3
\newdimen\DtwoY \DtwoY=\YshiftD \multiply\DtwoY by 3
\newdimen\EtwoX \EtwoX=\XshiftD \multiply\EtwoX by 4
\newdimen\EtwoY \EtwoY=\YshiftD \multiply\EtwoY by 4
\newdimen\FtwoX \FtwoX=\XshiftD \multiply\FtwoX by 5
\newdimen\FtwoY \FtwoY=\YshiftD \multiply\FtwoY by 5
\newdimen\BthreeX \BthreeX=0pt 
\newdimen\BthreeY \BthreeY=-\edgelengthD 
\newdimen\CthreeX \CthreeX=0pt
\newdimen\CthreeY \CthreeY=-\edgelengthD \multiply\CthreeY by 2
\newdimen\DthreeX \DthreeX=0pt
\newdimen\DthreeY \DthreeY=-\edgelengthD \multiply\DthreeY by 3
\newdimen\EthreeX \EthreeX=0pt
\newdimen\EthreeY \EthreeY=-\edgelengthD \multiply\EthreeY by 4
\newdimen\FthreeX \FthreeX=0pt
\newdimen\FthreeY \FthreeY=-\edgelengthD \multiply\FthreeY by 5
\bond\AX\AY\FoneX\FoneY
\bond\AX\AY\FtwoX\FtwoY
\bond\AX\AY\FthreeX\FthreeY
\hollownode\AX\AY
\hollownode\BoneX\BoneY
\hollownode\ConeX\ConeY
\hollownode\DoneX\DoneY
\hollownode\EoneX\EoneY
\hollownode\FoneX\FoneY
\hollownode\BtwoX\BtwoY
\hollownode\CtwoX\CtwoY
\hollownode\DtwoX\DtwoY
\hollownode\EtwoX\EtwoY
\hollownode\FtwoX\FtwoY
\hollownode\BthreeX\BthreeY
\hollownode\CthreeX\CthreeY
\hollownode\DthreeX\DthreeY
\hollownode\EthreeX\EthreeY
\hollownode\FthreeX\FthreeY
%
%\nearnode\AX\AY{0}{120}{0pt}{0pt}{b}{$a$}%
%\nearnode\BoneX\BoneY{-70}{-70}{0pt}{0pt}{tr}{$b_1$}%
%\nearnode\ConeX\ConeY{-70}{-70}{0pt}{0pt}{tr}{$c_1$}%
%\nearnode\DoneX\DoneY{-70}{-70}{0pt}{0pt}{tr}{$d_1$}%
%\nearnode\EoneX\EoneY{-70}{-70}{0pt}{0pt}{tr}{$e_1$}%
%\nearnode\FoneX\FoneY{-70}{-70}{0pt}{0pt}{tr}{$f_1$}%
%
%\nearnode\BtwoX\BtwoY{70}{-70}{0pt}{0pt}{tl}{$b_2$}%
%\nearnode\CtwoX\CtwoY{70}{-70}{0pt}{0pt}{tl}{$c_2$}%
%\nearnode\DtwoX\DtwoY{70}{-70}{0pt}{0pt}{tl}{$d_2$}%
%\nearnode\EtwoX\EtwoY{70}{-70}{0pt}{0pt}{tl}{$e_2$}%
%\nearnode\FtwoX\FtwoY{70}{-70}{0pt}{0pt}{tl}{$f_2$}%
%
%\nearnode\BthreeX\BthreeY{120}{0}{0pt}{0pt}{l}{$b_3$}%
%\nearnode\CthreeX\CthreeY{120}{0}{0pt}{0pt}{l}{$c_3$}%
%\nearnode\DthreeX\DthreeY{120}{0}{0pt}{0pt}{l}{$d_3$}%
%\nearnode\EthreeX\EthreeY{120}{0}{0pt}{0pt}{l}{$e_3$}%
%\nearnode\FthreeX\FthreeY{120}{0}{0pt}{0pt}{l}{$f_3$}%
%
\end{picture}
\end{displaymath}
That is, two generators braid ($aba=bab$) or commute ($ab=ba$) when
the corresponding nodes are joined or unjoined.  Basak \cite{basak} has
proven this, his proof making essential use of a computer.  Our
purpose is to give a conceptual, computer-free proof.  We hope that it
will clarify which structures will be important for further work
on the conjecture of \cite{monstrous}.

We will describe the players in the main theorem, state the theorem,
and then summarize the sections.  $\G$ is the group of all isometries
of a certain lattice $L_{13,1}$ over the Eisenstein integers
$\E=\Z[\w{=}\sqrt[3]{1}]$.  
This lattice has the property that
$L_{13,1}=\theta L_{13,1}'$, 
where $\theta=\w-\wbar=\sqrt{-3}$ and the
prime indicates the dual lattice.  Also, $L_{13,1}$ is the unique
$\E$-lattice of signature $(13,1)$ with this property.  An explicit
model for it appears in section~\ref{sec-L13,1}.  The Artin group of the
$Y_{555}$ diagram means the abstract group with one generator for each
node of the diagram, with the braid and commutation relations
described above.  Such groups arise naturally in the fundamental
groups of hyperplane complements.  A triflection means a complex
reflection of order 3, where a complex reflection means a nontrivial
isometry of a Hermitian vector space that fixes a hyperplane
pointwise.  Complex reflections arise naturally when studying branched
covers, in a manner explained in \cite{monstrous}.

\begin{theorem}[\cite{basak}]
\label{thm-Y555-generates}
Up to complex conjugation, there is a unique irreducible action of the
$Y_{555}$ Artin group on a Hermitian vector space of dimension${}>1$
in which the generators act by triflections.  The image of this
representation is $\aut L_{13,1}$.
\end{theorem}

In section~\ref{sec-E-lattices} we give background on Eisenstein lattices, and in
section~\ref{sec-root-and-Niemeier} we classify two types of such lattices, the root
lattices (analogous to the ADE lattices over $\Z$) and the Eisenstein
Niemeier lattices (equivalently, $\E$-lattice structures on the
Niemeier lattices).  The point of this is to enumerate the 5 cusps of
$\ch^{13}/P\G$ and be able to recognize one as having ``Leech type''.
Section~\ref{sec-L13,1} describes $L_{13,1}$ in a manner convenient for the
proof of the theorem, which appears in section~\ref{sec-generation-of-Aut(L13,1)}.  Throughout, we
use ATLAS notation \cite{atlas} for group extensions:  $A.B$, $A{:}B$ and
$A{\cdot}B$.  

\section{Eisenstein lattices}
\label{sec-E-lattices}

We have already introduced the Eisenstein integers $\E=\Z[\w]$ and
defined $\theta=\w-\wbar=\sqrt{-3}$.  
An $\E$-lattice $L$ means a free
$\E$-module equipped with an $\E$-valued Hermitian form $\ip{}{}$,
linear in its first argument and antilinear in its second.  The norm
$|x|^2$ of a vector means $\ip{x}{x}$.  We call $L$ nondegenerate if
$L^{\perp}=0$; in this case the dual lattice $L'$ means the set of all
$v\in L\tensor\C$ with $\ip{L}{v}\sset\E$.  All the lattices we will
meet satisfy $L\sset\theta L':=\theta\cdot(L')$, which is to say that
all inner products are divisible by $\theta$.  This should be thought
of as an ordinary integrality condition, because it means that the
underlying $\Z$-lattice $L^\Z$, with $x\cdot
y=\frac{2}{3}\Re\ip{x}{y}$, is integral and even.  The rescaling by
$\frac{2}{3}$ is not important; it is a nuisance arising from the fact that the
smallest scale at which $L$ is integral as an $\E$-lattice is
different from the smallest scale at which it is an integral
$\Z$-lattice.  Most of our lattices will also satisfy $L=\theta L'$,
which is the same as the unimodularity of $L^\Z$.  If $L=\theta L'$
then $\det L=\pm\theta^{\dim L}$ (which makes sense since $\dim L$
turns out to be even).

Examples of $L$ with $L=\theta L'$ are the Eisenstein versions of the
$E_8$ lattice (\cite[ch.~7, example~11b]{splag} or
theorem~\ref{thm-root-lattices} below) and the Leech lattice
(\cite{wilson}, scaled to have minimal norm~$6$).  It is
well-known that an indefinite even unimodular $\Z$-lattice is
determined by its signature, and there is a corresponding result for
$\E$-lattices.  Namely, an $\E$-lattice $L$ of signature $(p,n)$
satisfying $L=\theta L'$ exists if and only if $p-n\cong 0$
modulo~$4$, and $L$ is unique when this signature is indefinite.  A
proof appears in \cite{basak}.
The main player in this paper is this lattice of signature $(13,1)$,
for which we will write $L_{13,1}$.  We studied it in \cite{ch13},
using slightly different conventions (signature $(1,13)$ and $\ip{}{}$
linear in its second argument rather than its first) and a particular
explicit model.  In section~\ref{sec-L13,1} we will give a different
explicit model.

If $L$ is an $\E$-lattice satisfying $L\sset\theta L'$, then $r\in L$
is called a root of $L$ if $|r|^2=3$.  The language reflects two
things.  First, $r$ becomes a root in the usual sense (a vector of
norm~2) when we pass to $L^\Z$.  Second, the complex reflection
\begin{equation}
\label{eq-complex-reflection}
x\mapsto x+(\w-1)\frac{\ip{x}{r}}{|r|^2}r
\end{equation}
is an isometry of $L$, so that roots give reflections, analogously to
roots in $\Z$-lattices.  But this is a triflection; we call it the $\w$-reflection in $r$, since it
multiplies $r$ by $\w$ and fixes $r^{\perp}$ pointwise.  If $r$ and
$r'$ are nonproportional roots, then their $\w$-reflections braid
if and only if
$|\ip{r}{r'}|^2=3$.   One can check this by multiplying out $2\times2$
matrices.  

\section{Root lattices;  Niemeier lattices; Null vectors of
  $L_{13,1}$}
\label{sec-root-and-Niemeier}

At a key point in section~\ref{sec-generation-of-Aut(L13,1)} we will
need to recognize a particular null vector $\p$ of $L_{13,1}$ as
having ``Leech type'', which is to say that $\pperp/\spanof{\p}$ is a
copy of the complex Leech lattice.  The reader may skip this section
if he is prepared to accept one consequence of
theorem~\ref{thm-Neimeier-lattices} below: a primitive null vector of
$L_{13,1}$ whose stabilizer contains a copy of $L_3(3)$ has Leech
type.  There is a quicker-and-dirtier proof than the one we give, but
we think the $\E$-lattice classifications are interesting in themselves.

We will need to understand the orbits of primitive null
vectors in $L_{13,1}$.  These turn out to be in bijection with the
positive-definite 12-dimensional lattices $L$ satisfying $L=\theta
L'$; we will call such lattices Eisenstein Niemeier lattices, since
their real forms are positive-definite 24-dimensional even unimodular
lattices, classified by Niemeier.  Since root lattices play a major
role in Niemeier's classification, they do in ours too, so we define
an Eisenstein root lattice to be a positive-definite $\E$-lattice $L$
satisfying $L\sset\theta L'$ and spanned by its roots.  

We will establish the bijection between Eisenstein
Niemeier lattices and orbits of primitive null vectors in $L_{13,1}$,
then classify the Eisenstein root lattices, and then use this to
classify the Eisenstein Niemeier lattices.  The root lattice
classification is similar to and simpler then the well-known ADE
classification of root lattices over $\Z$.  The Eisenstein Niemeier
lattices turn out to correspond to five of the classical Niemeier
lattices.

\begin{lemma}
\label{lem-cusps-and-null-vecs}
Suppose $p,n>0$, $p-n\cong 0\ (4)$, and $L_{p,n}$ is the unique
$\E$-lattice of signature $(p,n)$ satisfying $L_{p,n}=\theta
L_{p,n}'$.  If $\p$ is a primitive null vector then $L:=\p^{\perp}/\spanof{\p}$ is a lattice of signature $(p-1,n-1)$ that satisfies
$L=\theta L'$.  Every such $L$ arises this way.  Two primitive null
vectors $\p_1,\p_2$ of $L_{p,n}$ are equivalent under $\aut L_{p,n}$ if
and only if $\p_1^{\perp}/\spanof{\p_1}\iso \p_2^{\perp}/\spanof{\p_2}$.
\end{lemma}

\begin{proof}
(This is essentially the same as for even unimodular $\Z$-lattices.)
  By $L_{p,n}=\theta L_{p,n}'$, there exists $w\in L$ with
  $\ip{\p}{w}=\theta$.  Adding a multiple of $\p$ to $w$ allows us to
  also assume $|w|^2=0$, so $\spanof{\p,w}\iso\bigl(\begin{smallmatrix}0&\theta\\\thetabar&0\\\end{smallmatrix}\bigr)$.
  Therefore $\spanof{\p,w}=\theta\spanof{\p,w}'$, so
  $\spanof{\p,w}$ is a summand of $L_{p,n}$.  The other summand
  $\spanof{\p,w}^{\perp}$ must also satisfy $\spanof{\p,w}^{\perp}=\theta\bigl(\spanof{\p,w}^{\perp}\bigr)'$, and
  it projects isometrically to $\p^{\perp}/\spanof{\p}$.  This
  establishes the first claim.  For the second, given $L$ of signature
  $(p-1,n-1)$ satisfying $L=\theta L'$, we have
  $L\oplus \bigl(\begin{smallmatrix}0&\theta\\\thetabar&0\\\end{smallmatrix}\bigr)\iso
    L_{p,n}$, and it is now obvious that $L$ is $\p^{\perp}/\spanof{\p}$ for a suitable null vector $\p$.  In the last claim, if
    $\p_1$ and $\p_2$ are equivalent, then obviously $\p_1^{\perp}/\spanof{\p_1}\iso \p_2^{\perp}/\spanof{\p_2}$, so it suffices to
    show the converse.  The argument for the first claim implies that
    there is a direct sum decomposition $L_{p,n}\iso
    M_1\oplus\bigl(\begin{smallmatrix}0&\theta\\\thetabar&0\\\end{smallmatrix}\bigr)$ with
    $M_1\iso \p_1^{\perp}/\spanof{\p_1}$ and $\p_1$ corresponding to
    one of the coordinate vectors of the $2\times 2$ block.  And
    there is a similar decomposition with $\p_2$ in place of $\p_1$.
    Then, given an isomorphism $M_1\iso M_2$, it is easy to write down
    an automorphism of $L_{p,n}$ sending $\p_1$ to $\p_2$.
\end{proof}

\begin{theorem}
\label{thm-root-lattices}
Any Eisenstein root lattice is a direct sum of copies of the following
4 lattices:
\begin{align*}
A_2^\E&{}=\theta\E\\
D_4^\E&{}=\bigl\{(x,x,y)\in\E^3:x\cong y\ (\theta)\bigr\}\\
E_6^\E&{}=\bigl\{(x,y,z)\in\E^3:x\cong y\cong z\ (\theta)\bigr\}\\
E_8^\E&{}=\bigl\{(x_1,\dots,x_4)\in\E^4:\pi(x_1,\dots,x_4)\in\code_4\sset\F_3^4\bigr\},
\end{align*}
which have the properties listed in table~\ref{tab-root-lattices}.
(We use the standard inner product on $\C^n$.  Also, the description
of $E_8^\E$ refers to the map $\pi:\E^4\to\E^4/\theta\E^4=\F_3^4$ and
the tetracode $\code_4$, i.e., the subspace of $\F_3^4$ spanned by
$(0,1,1,1)$ and $(1,0,1,-1)$.)
\end{theorem}

\begin{table}
\label{tab-root-lattices}
\begin{tabular}{cccrcc}
&&&&&$\theta L'-L$\\
\medskip % mysteriously, puts this skip AFTER the following line
$L$&$R$&$\aut L$&$|\aut L|$&
$\theta L'/L$&min. norm\\
$A_2^\E$&$\Z/3$&${}\times\Z/2$&6&$\F_3^1$&1\\
$D_4^\E$&$\SL_2(3)$&${}\times\Z/3$&72&$\F_4^1$&3/2\\
$E_6^\E$&$3^{1+2}{:}\SL_2(3)$&${}\times\Z/2$&1,296&$\F_3^1$&2\\
$E_8^\E$&$3\times\Sp_4(3)$&${}\times1$&155,520&$0$&\\
\end{tabular}
\bigskip
\caption{The indecomposable Eisenstein root lattices.  The second
  column gives the structure of the group $R$ generated by the
  triflections in the roots of $L$, in ATLAS notation \cite{atlas}.
  $\aut L$ is the product of this group with the cyclic group of
  scalars given in the third column.  The fifth column describes
  $\theta L'/L$ as a vector space over $\E/\theta\E=\F_3$ or
  $\E/2\E=\F_4$.  Every nonzero element of $\theta L'/L$ has minimal
  representatives of norm given in the last column.}
\end{table}

\begin{proof}
The data in the table will be helpful in the classification, so we
begin there.  The claims for $L=A_2^\E$ are obvious; we remark that
the smallest elements of $\theta L'-L$ are the units of $\E$, and all
others have norm${}>3$.  

Now let $L=D_4^\E$.  Its 24 roots are the scalar multiples of
$(\w^i,\w^i,1)$ and $(0,0,\theta)$.  It is easy to see that
conjugation by each of
the $24/6=4$ cyclic groups generated by triflections permutes the
other 3 cyclically.  Therefore $R$ is generated by two triflections
that braid, so it is an image of $\langle
a,b\,|\,aba=bab,\ a^3=b^3=1\rangle$, which is a presentation for
$\SL_2(3)$.  To see that $R$ is $\SL_2(3)$ rather than a proper
quotient, consider its action on $L/\theta L\iso\F_3^2$.  Now, $R$
permutes the scalar classes of roots as the alternating group $A_4$,
so if we choose any 2 non-proportional roots $r$ and $s$, then $\aut L$
is generated by $R$ together with the transformations sending $r$ to a
multiple of itself and $s$ to a multiple of itself.  Since
$\ip{r}{s}\neq0$, such a transformation must be a scalar.  So $\aut
L=R\times\spanof{\w}$.   Finally, it is easy to see that the norm
6 vectors of $L$ are the scalar multiples of $(\w^i,\w^i,1+\theta)$
and $(\theta,\theta,0)$, and that the halves of these vectors span
$\theta L'$.  In fact, the halves of these vectors account for all the
elements of $\theta L'-L$ of norm${}\leq3$, and are all equivalent under
$\aut L$.  Representatives for $\theta L'/L$ are $0$ and
$\frac{\w^i}{2}(\theta,\theta,0)$, so $\theta L'/L\iso\F_4^1$.  

Now let $L=E_8^\E$.  Because it is got from $(A_2^\E)^4$ by
gluing along the 2-dimensional code
$\code_4\sset\bigl(\theta(A_2^\E)'/A_2^\E\bigr)^4\iso\F_3^4$, it
satisfies $L=\theta L'$, justifying the last two entries in the
table.  The descriptions of $R$ and $\aut L$ are justified by
theorem~5.2 of \cite{ch13}.  (The proof in \cite{ch13} appeals to a
coset enumeration to establish that $L$ contains the scalars of
order~3; this may be avoided by observing that $L$ contains 4 mutually
orthogonal roots.)

Now let $L=E_6^\E$ and note that the following symmetries are visible:
permutation of coordinates, multiplication of coordinates by cube
roots of unity, and the scalar~$-1$.  It is easy to see that $\theta
L'=\{(x,y,z)\in\E^3:x+y+z\cong0\ (\theta)$\}, whose 54 minimal vectors
are got from $(1,-1,0)$ by applying these symmetries.  Note also that
these are the only elements of $\theta L'-L$ of norm${}\leq3$.  It is
easy to see that $\theta L'/L\iso\F_3^1$.  Also, $L$ is the orthogonal
complement of $r=(0,0,0,\theta)\in E_8^\E$, and every automorphism
$\phi$ of $E_6^\E$ extends uniquely to an automorphism of $E_8^\E$
that either fixes or negates $r$.  (The extension fixes or negates $r$
according to whether $\phi$ fixes or negates $\theta
L'/L\iso\theta\spanof{r}'/\spanof{r}\iso\F_3^1$.)  It
follows that $|\aut L|\leq 2\times\frac{1}{240}\times 155,520=1296$.
Since triflections must act trivially on $\F_3^1$, we also have
$|R|\leq648$.  We will show that $R$ has structure
$3^{1+2}{:}\SL_2(3)$;
this will justify the first column of the table, and (since $-1\notin
R$) also the second.  

To see the map $R\to\SL_2(3)$, consider the action on
$L/3L'\iso\F_3^2$.  All roots are equivalent under $\aut L$ (since any
two root in a $D_4^\E$ are $R$-equivalent), and the 72 roots fall into
8 classes of size 9, accounting for all 8 nonzero elements of $L/3L'$.
This space supports a symplectic form, given by dividing inner
products by $\theta$ and then reducing mod $\theta$.  The
$\w$-reflection in a root projects to the symplectic transvection in
the image of the root.  Now we study the kernel $K$ of $R\to\SL_2(3)$.
If $r$ and $s$ are orthogonal roots then their $\w$-reflections map to
the same transvection $L/3L'$ (since they map to commuting
transvections), so the quotient of the reflections lies in $K$.  This
shows: if an automorphism of $L$ has 3 roots as eigenvectors, with
eigenvalues $1$, $\w$ and $\wbar$, then it lies in $K$.  For example,
$\diag[1,w,\wbar]\in K$.  Also, the cyclic permutation of coordinates
lies in $K$.  These two elements of $K$ generate an extraspecial group
$3^{1+2}$.  Since $27\cdot|\SL_2(3)|=648$, we have shown
$R=3^{1+2}.\SL_2(3)$.  The extension splits because the reflection
group of a $D_4^\E$ sublattice provides a complement.

Having established the table, we will now classify the Eisenstein root
lattices.  Call such a lattice decomposable if its roots fall into two
or more nonempty classes, with members of distinct classes being orthogonal.
In this case it is a direct sum of lower-dimensional root lattices, so
it suffices to show that $A_2^\E$, $D_4^\E$, $E_6^\E$ and $E_8^\E$ are
the only indecomposable Eisenstein root lattices.  We will use the
following facts, established above.  (i)~If $L=A_2^\E$, $D_4^\E$ or
$E_6^\E$, then $\aut L$ acts transitively on the vectors of $\theta
L'-L$ of norm${}\leq3$. (ii)~$E_8^\E=\theta(E_8^\E)'$.  

Suppose $M$ is an indecomposable Eisenstein root lattice.  If $\dim M=1$ then
obviously $M\iso A_2^\E$.  If $\dim M=2$ then it contains a
1-dimensional indecomposable root lattice $L$, and we know $L\iso
A_2^\E$.  Also, $M$ contains a root $r$ not in $L\tensor\C$, whose
projection to $L\tensor\C$ is nonzero.  Since this projection is an
element $r$ of $\theta L'-\{0\}$ of norm${}<3$, and $\aut L$ acts
transitively on such vectors, there is an essentially unique
possibility for $\spanof{L,r}$.  Since $D_4^\E$ arises by this
construction, $\spanof{L,r}\iso D_4^\E$.  Therefore $M$ lies
between $\theta(D_4^\E)'$ and $D_4^\E$.  Since every norm 3 vector of
$\theta(D_4^\E)'$ lies in $D_4^\E$, we have $M\iso D_4^\E$.   If $\dim
M=3$ then the same argument, with $L=D_4^\E$, shows that $M\iso
E_6^\E$.  If $\dim M>3$, then the same argument, with $L=E_6^E$, shows
that $M$ contains $E_8^\E$.  Then $E_8=\theta(E_8^\E)'$ implies that
$E_8^\E$ is a summand of $M$, and indecomposability implies $M\iso
E_8^\E$. 
\end{proof}

\begin{theorem}
\label{thm-Neimeier-lattices}
There are exactly 5 Eisenstein Niemeier lattices:

$\bigl(A_2^\E\bigr)^{12}$ glued along the ternary Golay code,
with group $3^{12}{:}2M_{12};$

$\bigl(D_4^\E\bigr)^6$ glued along the hexacode, with group $\SL_2(3)^6{:}3A_6;$

$\bigl(E_6^\E\bigr)^4$ glued along the tetracode, with group $\bigl(3^{1+2}{:}\SL_2(3)\bigr)^4{:}\SL_2(3);$

$\bigl(E_8^\E\bigr)^3$, with group
$\bigl(3\times\Sp_4(3)\bigr)^3{:}S_3;$ and

the complex Leech lattice $\cll$, with group $6\suz$.
\end{theorem}

\noindent 
Here, $M_{12}$ and $\suz$ are the sporadic finite simple groups of
Mathieu and Suzuki.

\begin{proof}
Our argument is similar in spirit to Venkov's treatment \cite{venkov}
of Niemeier's classification.  Suppose $L$ is an Eisenstein Niemeier
lattice and $L^\Z$ its underlying real lattice.  By Niemeier's
classification, there are 24 possibilities for $L^\Z$; in 23 cases the
roots span $L^\Z$ up to finite index, and in the last case $L^\Z$ has
no roots and is the Leech lattice.  By
theorem~\ref{thm-root-lattices}, the root system of $L^\Z$ must be a
sum of $A_2$, $D_4$, $E_6$ and $E_8$ root systems.  Considering
Niemeier's list shows that $L^\Z$'s root system is $A_2^{12}$,
$D_4^6$, $E_6^4$, $E_8^3$ or empty.  We treat the first four cases
first.  Theorem~\ref{thm-root-lattices} shows that there is a unique
Eisenstein structure on the root sublattice of $L^\Z$, so the
sublattice $L_0$ of $L$ spanned by its roots is $(A_2^\E)^{12}$,
$(D_4^\E)^6$, $(E_6^\E)^4$ or $(E_8^\E)^3$.  In the last case we have
$L=L_0$ and are done.  In the other cases, $L$ lies between $\theta
L_0'$ and $L_0$, so it is determined by its image $\code$ in $\theta
L_0'/L_0\iso\F_3^{12}$, $\F_4^6$ or $\F_3^4$ in the three cases.  We
must have $\code\sset\codeperp$ (with respect to the usual quadratic
form on $\F_3^n$ or Hermitian form on $\F_4^6$), in order to have
$L\sset\theta L'$.  Also, $\code$ must be half-dimensional in $\theta
L_0'/L_0$, in order to have $L=\theta L'$.  Finally, all roots of $L$
already lie in $L_0$, by definition.

In the $A_2$ case, these conditions imply that $\code$ is a selfdual
code of length 12 with no codewords of weight~3.  The ternary Golay
code is the unique such code, up to monomial transformations of
$\F_3^{12}$, so $\code$ is it and $L$ is as described.  In the $D_4$
case, $\code$ is a selfdual subspace of $\F_4^6$ with no codewords of
weight~2.  The hexacode is the unique such code, up to monomial
transformations, so $\code$ is it and $L$ is as described.  In the
$E_6$ case, $\code$ is a 2-dimensional subspace of $\F_3^4$ having no
codewords of weight${}<3$.  Again there is a unique candidate, the
tetracode, and $L$ is as described.

Next we treat the case that $L^\Z$ is the Leech lattice; we must show
that $L$ is the complex Leech lattice.  I know of 3 completely
independent approaches.  (1) The uniqueness of the $\E$-module
structure on the Leech lattice is the same as the uniqueness of the
conjugacy class in $\Co_0=\aut(L^\Z)$ of elements of order 3 with no
fixed vectors.  This can be checked by consulting the character table
\cite{atlas} for $\Co_0$.  (2) Use
lemma~\ref{lem-cusps-and-null-vecs}, together with theorem~4.1 of
\cite{ch13}, which contains the statement that $L_{13,1}$ has a unique
orbit of primitive null vectors orthogonal to no roots.  (3) Presumably
one can mimic Conway's characterization of the Leech lattice
\cite{conway-leech-characterization}, applying analogues of his
counting argument to $L/\theta L$.

The automorphism group of $\cll$ is treated in detail in
\cite{wilson}.  The other automorphism groups are easy to work out.
Let $L_0=M^n$ be the decomposition of $L_0$ into its indecomposable
summands and let $R$ be the group generated by triflections in the roots
of $M$.  Recall from theorem~\ref{thm-root-lattices} that $\aut M$ splits as
$R\times C$, where $C$ denotes the group of scalars from column~3 of
table~\ref{tab-root-lattices}.  Obviously $\aut L\sset\aut L_0=(R^n\times C^n){:}S_n$;
indeed it is the subgroup of this that preserves $\code\sset(\theta
M'/M)^n$.  Now, $R$ acts trivially and $C^n{:}S_n$ acts by monomial
transformations.  Therefore $\aut L$ is the semidirect product of
$R^n$ by the subgroup of $C^n{:}S_n$ whose action preserves $\code$.
This latter group is $2M_{12}$, $3A_6$,
$\SL_2(3)$ or $S_3$ in the four cases.  (The
automorphism group of the hexacode is sometimes given as
$3{\cdot}S_6$, but the elements not in $3A_6$ are $\F_4$-antilinear,
so they arise from antilinear maps $L\to L$.)
\end{proof}

\section{A model of $L_{13,1}$}
\label{sec-L13,1}

In this section we describe $L_{13,1}$ with $3^{13}{:}L_3(3)$ among
its visible symmetries.  We begin with the vector space
$\F_3^{13}$, with coordinates indexed by the points of $P^2\F_3$, and
proceed to define two codes.  The first is the ``line difference
code'' $\code$, spanned by the differences of (characteristic
functions of) lines of $P^2\F_3$, and the second is the ``line code'',
which derives its name from the fact that it is spanned by lines, but
is formally defined (and written) as $\codeperp$ (with respect to the
usual inner product).

Two lines of $P^2\F_3$ meet in 1 point (or 4), and it follows that
$\code$ is orthogonal to every line, hence orthogonal to itself.
Therefore $\dim\code\leq6$.  On the other hand, it is easy to
enumerate some elements of $\code$ (table~\ref{tab-elements-of-C}).
This shows that $\dim\code=6$ and also that the enumeration is complete.
Therefore $\dim\codeperp=7$, and since a line lies in $\codeperp$ but
not $\code$, we see that $\codeperp$ is indeed spanned by lines.  It
will be useful to have a list of the elements of $\codeperp$: these are
the codewords in table~\ref{tab-elements-of-Cperp}, their negatives,
and the elements of $\code$.  We compiled
table~\ref{tab-elements-of-Cperp} by adding the all 1's vector (the
sum of all 13 lines) to the elements of $\code$.

\begin{table}
\begin{tabular}{cclr}
\medskip
support&coordinates&description&\llap{number}\\
0&$0^{13}$&&1\\
6&${+}^3{-}^30^7$&difference of lines&156\\
9&$\pm({+}^90^4)$&affine plane&26\\
9&$\pm({+}^6-^30^4)$&sum of 3 general lines&468\\
12&${+}^6-^60$&$l_1+l_2-l_3-l_3$ for 4 concurrent lines&78\\
\end{tabular}
\bigskip
\caption{The elements of $\code$.}
\label{tab-elements-of-C}
\end{table}

\begin{table}
\begin{tabular}{cclr}
\medskip
support&coordinates&description&\llap{number}\\
4&${+}^40^9$&line&13\\
7&${+}^1-^60^6$&sum of two lines (negated)&78\\
7&${+}^4-^30^6$&$-1$ on vertices, $0$ on edges, $1$ elsewhere&234\\
10&${+}^4-^60^3$&$0$ on vertices, $-1$ on edges, $1$ elsewhere&234\\
10&${+}^7-^30^3$&$-1$ on $l_1-l_2$, $0$ on $l_2-l_1$, $1$ elsewhere&156\\
13&${+}^4-^9$&$1$ on a line, $-1$ elsewhere&13\\
13&${+}^{13}$&&1\\
\end{tabular}
\bigskip
\caption{The elements of $\codeperp$ with coordinate sum~$1$. 
The entries that
refer to ``vertices'' and ``edges'' refer to three general lines---a
vertex means a point on two of the lines, and an edge means a point on
just one of them.}
\label{tab-elements-of-Cperp}
\end{table}

We work with the usual inner product of signature $(13,1)$ on $\C^{14}$,
$$
\ip{x}{y}=-x_0\bar{y}_0+x_1\bar{y}_1+\dots+x_{13}\bar{y}_{13},
$$
index the last 13 coordinates by the points of $P^2\F_3$, and define
$L$ as the set of vectors $(x_0;x_1,\dots,x_{13})$ such that 
$x_0\cong x_1+\dots+x_{13}$ mod $\theta$ and that $(x_1,\dots,x_{13})$,
modulo $\theta$, is an element of $\codeperp$.

\begin{theorem}
\label{thm-model-of-L13,1}
$L$ is isomorphic to $L_{13,1}$ and is spanned by the $13$ ``point
roots'' $(0;\theta,0^{12})$, with the $\theta$ in any of the last $13$
spots, and the $13$ ``line roots'' $(1;1^4,0^9)$, with the $1$'s
along a line of $P^2\F_3$.  
\end{theorem}

\begin{proof}
It is easy to see that the point and line roots span $L$.  If $p$ is a
point root and $\ell$ a line root, then $\ip{p}{\ell}=\theta$ or $0$
according to whether the point lies on the line.  Also, any two point
roots are orthogonal, as are any two line roots.  Therefore
$L\sset\theta L'$.  To see $L=\theta L'$, check that $L$ contains
$(\theta;0^{13})$ and consider the span $M$ of it and the point roots.
Then
$\theta M'/M\iso\F_3^{14}$, and we need to check that the image of $L$
therein is $7$-dimensional.  This is easy because we know
$\dim\codeperp=7$ and the $0$th coordinate of an element of $L$
is determined modulo $\theta$ by the others.
\end{proof}

The promised group $3^{13}{:}L_3(3)$ is generated by the triflections
in the point roots and the permutations of the last 13 coordinates by
$L_3(3)$.  

The following lemma is not central; it is used only to
establish the equality of two lattices in the proof of
lemma~\ref{lem-specific-roots-in-R}.

\begin{lemma}
\label{lem-unique-enlargement}
Let $M$ be the $12$-dimensional lattice consisting of all vectors in
$(\theta\E)^{13}$ with coordinate sum zero.  Then there is a
unique lattice $N$ preserved by $L_3(3)$, strictly containing $M$, and
satisfying $N\sset\theta N'$.
\end{lemma}

\begin{proof}[Proof sketch:\/]
  Any lattice $N$ containing $M$ and satisfying $N\sset\theta N'$ lies
  in $\theta M'$, so that it corresponds to a subspace of $Z:=\theta
  M'/M$.  And $Z$ is the coordinate-sum-zero subspace of $\F_3^{13}$.   The
  lemma follows from the fact that $\code$ is the unique
  nontrivial $L_3(3)$-invariant subspace.  To see this, one checks
  that $\code$ is irreducible under $L_3(3)$, so that $Z/\code$
  is also irreducible (being the dual), and that
  $\code$ has no invariant complement.
\end{proof}

\section{Generation of $\aut(L_{13,1})$ by the $Y_{555}$ triflections}
\label{sec-generation-of-Aut(L13,1)}

In this section we prove the main theorem, theorem~\ref{thm-Y555-generates}.  First we
prove uniqueness.  Label the generators by $g_1,\dots,g_{16}$.  The
argument of \cite[sec.~5]{ch13} shows that without loss we may take
the $g_i$ to be the $\w$-reflections in pairwise linearly
independent vectors $r_1,\dots,r_{16}$ of norm $3$, satisfying
$|\ip{r_i}{r_j}|^2=3$ or $0$ according to whether $g_i$ and $g_j$
braid or commute.  It is convenient to 2-color $Y_{555}$ and suppose
$\ip{r_i}{r_j}=\theta$ (resp. $-\theta$) when $g_i$ and $g_j$ braid
and $g_i$ is black (resp. white).  The inner product matrix of the
$r_i$ 
turns out to have rank~14 (by
direct computation or the realization below), so $V$ must have
dimension~14 (by irreducibility of $V$ and connectedness of $Y_{555}$).  The $r_i$ are determined up to
isometries of $V$ by their inner products, so their configuration is
unique.

Having proven uniqueness of the representation, we can define $R$ as
its image.  To identify $R$ with $\aut L_{13,1}$, we will use the
model for $L_{13,1}$ from the previous section, and write $L$ for it.
Let $\D$ be the incidence graph of the points and lines of $P^2\F_3$,
and color the nodes corresponding to points black
and lines white.  Then the point and line roots from theorem~\ref{thm-model-of-L13,1}
satisfy the same inner product conditions as the $r_i$ chosen above.
It is possible (uniquely up to $L_3(3)$) to embed the $Y_{555}$
diagram into $\D$, preserving node colors.  So the 16 roots for
$Y_{555}$ may be taken to be 16 of the point and line roots.  It would
be annoying to make a choice of which 16, and we are saved from this
by the following lemma.

\begin{lemma}
\label{lem-Y555-roots-span-and-R-has-L3(3)}
The $16$ roots for $Y_{555}$ span $L$, and $R$ contains $L_3(3)$ and
the triflections in all the point and line roots.
\end{lemma}

\begin{proof}
First observe that $Y_{555}$ contains an 11-chain $E$ and a 4-chain
$F$ not joined to it.  By \cite[fig.~5.1]{ch13}, the roots of $E$ span a copy of $L_{9,1}$ and
those of $F$ a copy of $E_8^\E$, so
together they span $L$.  This proves our first assertion.  

One can check that for any 11-chain $E$ in $\D$, $E$ has a unique
extension to a 12-cycle $C$, and that the nodes of $\D$ not joined to
$C$ form a 4-chain $F$.  ($E$ is unique up to $L_3(3){:}2$, so
checking a single example suffices.)
We claim that if $R$ contains the triflections in the roots of $E$,
then it also contains the triflections in the root extending $E$ to
$C$.  We use a computation-free variation of the proof of
\cite[lemma~3.2]{basak}.  First use the fact that the roots of $F$ span a
copy of $E_8^\E$, whose orthogonal complement in $L$ must be a copy of
$L_{9,1}$.  By \cite[thm.~5.2]{ch13}, $\aut L_{9,1}$ is generated by
the triflections of $E$ and hence lies in $R$.  And since the
extending root is orthogonal to $E_8^\E$, it also lies in $L_{9,1}$, so its
triflections also lie in $R$.  This proves the claim.  Now, starting
with the three 11-chains in $Y_{555}$ and repeatedly applying the
claim shows that $R$ contains the triflections in all 26 roots.

We use a similar trick to show $L_3(3)\sset R$. Consider any
$Y_{555}\sset\D$, and let $E$ be one of its 11-chains and $F$ the
4-chain in $Y_{555}$ not joined to it.  Let $\phi$ be the diagram
automorphism of $Y_{555}$ that fixes each node of $F$ and exchanges
the ends of $E$.  One can check that $\phi$ extends to an automorphism
of $\D$, preserving node colors since it has a fixed point.  Therefore
$\phi$ defines an automorphism of $L$, permuting the point and line
roots as it permutes the points and lines of $P^2\F_3$.  Since $\phi$
fixes $F$ pointwise, it is an automorphism of the $L_{9,1}$ spanned by
the roots of $E$.  We already know that $R$ contains $\aut L_{9,1}$,
so it contains $\phi$.  So each $Y_{555}\sset\D$ gives rise to an
$S_3\sset R\cap L_3(3)$.  The set of elements of $R\cap L_3(3)$ so
obtained, from all $Y_{555}$ subdiagrams, is clearly normal in
$L_3(3)$.  Since $L_3(3)$ is simple, $R\cap L_3(3)$ is all of
$L_3(3)$.  (This diagram-automorphism trick was also used in
\cite[thm.~5.1]{ch13} and \cite[thm.~5.8]{basak}.)
\end{proof}

The next lemma shows that if $R$ contains certain triflections, then
it contains a well-understood group, of finite index in the stabilizer
of a null vector.  The lemma after that shows that $R$ does indeed
contain these triflections, and then we can complete the proof of
theorem~\ref{thm-Y555-generates} by showing $R=\aut L$.

\begin{lemma}
\label{lem-stabilizer-of-null-vector}
Suppose $L$ is an $\E$-lattice of dimension${}>2$ satisfying $L=\theta
L'$.  Suppose $\p$ is a primitive null vector and $r_i$ are
roots satisfying $\ip{r_i}{\p}=\theta$, such that the span of their
differences projects onto $\pperp/\spanof{\p}$.  Let $G$ be
the group generated by the triflections in the $r_i$ and $\p+r_i$.
Then $G$ contains every element of $\aut L$ that acts by a scalar on
$\spanof{\p}$ and trivially on $\pperp/\spanof\p$ 
\end{lemma}

\begin{remark}
The hypothesis $\dim L>2$ is necessary and should also have been imposed in theorem~3.2
of \cite{ch13}.
\end{remark}

\begin{proof}
This is implicit in the proofs of theorem~3.1 and~3.2 of \cite{ch13};
since the argument is slightly different and our conventions there
were different, we phrase the argument in coordinate-free language and
refer to \cite{ch13} for the supporting calculations.  By the
unipotent radical $U$ of the stabilizer of $\p$ we mean the automorphisms
of $L$ that fix $\p$ and act trivially on $M:=\pperp/\spanof\p$.  It
is a Heisenberg group, with center $Z$ equal to its commutator
subgroup and isomorphic to $\Z$, with $U/Z$ a copy of the additive
group of $M$.  The set $X$
of scalar classes of roots $r$ with $\bigl|\ip{r}{\p}\bigr|=|\theta|$
is a principal homogeneous space for $U$, and the set $X/Z$ of its
$Z$-orbits is a principal homogeneous space for $U/Z\iso M$.  If $r$
is a root with $\ip{r}{\p}=\theta$, then the triflections in $r$ and
$\p+r$ can be composed to yield a transformation multiplying $\p$ by a
primitive 6th root of unity and acting on $X/Z$ by scalar
multiplication by a primitive 6th root of unity, where $X/Z$ is
identified with $M$ by taking $rZ$ as the origin.   Write $\phi_r$
for this transformation  (which depends only on $rZ$, though we
don't need this).  

Suppose $r'$ is another root with $\ip{r'}{\p}=\theta$.  Since $X/Z$
is a principal homogeneous space modeled on $M$, there exists $m\in M$
with $m\cdot rZ=r'Z$.  Then $\phi_r\circ\phi_{r'}^{-1}$ turns out to
be an element of $U$, acting on $X/Z$ by translation by a unit times
$m$.  Under the hypothesis of the lemma, $G$ contains elements of $U$
for sufficiently many $m$ to span $M$ as an $\E$-lattice.  Taking
conjugates by (any) $\phi_r$ gives the unit scalar multiples of these
$m$, so $G$ contains enough elements of $U$ to generate $M$ as a
group.  Taking commutators shows that $G$ contains $Z$, so it contains
all of $U$.  And $\spanof{U,\phi_r}$ consists of all the elements of
$\aut L$ that we are asserting to lie in $G$.
\end{proof}

\begin{lemma}
\label{lem-specific-roots-in-R}
Let $\p$ be the primitive null vector $(-4-\w;1^{13})\in L$.  If $r$
is one of the $156$ roots $(2+\theta;0^3,\wbar^3,-1^7)$ or one of the $234$
roots $(-2\wbar;\wbar^4,-1^3,0^6)$, then $\ip{r}{\p}=\theta$ and $R$
contains the triflections in $r$ and $\p+r$.  The differences of these $390$
roots span $\pperp$.
\end{lemma}

\begin{remarks}
The exact placement of the coordinates can be determined up to
$L_3(3)$ by reducing the last 13 coordinates modulo $\theta$ and
comparing with the list of elements of $\codeperp$.  For example, for
one of the 156 roots, the $0$'s lie on one line of $P^2\F_3$, the
$\wbar$'s lie on another, and the $-1$'s are everywhere else,
including the point where the lines intersect.  (There are
$13\cdot12=156$ ways to choose the two lines.)  The same method
applies to all vectors referred to in the proof.  Also, $R$ contains
the triflections in some less-complicated roots $r$ satisfying
$\ip{r}{\p}=\theta$, for example the point roots.  But for these,
showing that $R$ contains the triflections in $\p+r$ is harder.  We
chose these roots because both $r$ and $\p+r$  have
small $0$th coordinate.
\end{remarks}

\begin{proof}
Checking $\ip{r}{\p}=\theta$ is just a computation.  Now we show that
$R$ has various roots $r$ (meaning that it contains the triflections
in them).  We will use the following fact repeatedly: if $R$ has roots
$a$ and $b$, and $\ip{a}{b}=\w-1$ or $\wbar-1$, then $a+b$ is a root
and $R$ has it too.  (This is because $\spanof{a,b}\iso D_4^\E$, and
the reflections in any two independent roots of $D_4^\E$ generate the
whole reflection group of $D_4^\E$.)  We know already that $R$ has the
line roots and their images under scalars and $3^{13}{:}L_3(3)$.

Step 1:  {\it $R$ has the roots $(\theta;1^3,-1^3,0^7)$ with the $1$'s
collinear and the $-1$'s collinear.}  Take $b$ to be the line root
$(1;1^4,0^9)$, and try $a$ having the form $(-\w;?,0^3,?^3,0^6)$,
where the $?$'s are negated cube roots of $1$, lying along a different
line.  We try this $a$ because $a+b=(1-\w;\dots)$, so if we can
choose the $?$'s with $\ip{a}{b}=\w^{\pm1}-1$, then we can conclude
that $R$ has a root $a+b=(1-\w;\dots)$, which we didn't know
before.  
We may in fact achieve $\ip{a}{b}=\w-1$, by taking (say) all the $?$'s to
be $-1$.  Then $R$ has the root $(1-\w;0,1^3,-1^3,0^6)$.  Applying a
scalar and an element of $3^{13}{:}L_3(3)$ finishes step~1 (this part
of the argument will be left implicit in steps~2--5).

Step 2: {\it $R$ has the roots $(2;-1^4,1^3,0^6)$ with the $1$'s at
  three noncollinear points, the $0$'s on the lines joining them, and
  the $-1$'s everywhere else.}  Take $b=(\theta;1^3,-1^3,0^7)$ from
step~1, and try $a=(1;?,0^2,?,0^2,?^2,0^5)$, with the $?$'s lying on a
line that meets a $1$ and a $-1$ of $b$.  Solving for the $?$'s as
before reveals that $a=(1;\wbar,0^2,1,0^2,?^2,0^5)$ satisfies
$\ip{a}{b}=\w-1$.  So $R$ has the root
$(\theta+1;-\w,1^2,0,-1^2,?^2,0^5)$ where the $?$'s are cube roots of
$1$---exactly which ones is unimportant.

Step 3: {\it $R$ has the roots $(2;1^6,-1,0^6)$, where the $1$'s all
  lie on two lines through the $-1$.}  Take $b=(1;1^4,0^9)$ and try
$a=(1;?,0^3,?^3,0^6)$.  Solving for the $?$'s reveals that
$a=(1;\w,0^3,?^3,0^6)$ satisfies $\ip{a}{b}=\w-1$.  So $R$ has the
root $a+b=(2;-\wbar,1^3,?^3,0^6)$ with the $?$'s being cube roots
of~$1$.

Step 4: {\it $R$ has the roots $(2-\wbar;-1^3,0^3,1^7)$, with the
  $-1$'s collinear and the $0$'s collinear.}  Take $b=(2;1^6,-1,0^6)$
from step~3, and try $a=(-\wbar;0^6,?,?^3,0^3)$ where the $?$'s all
lie on a line through the $-1$ of $b$.  Solving for the $?$'s reveals
that $a=(-\wbar;0^6,-\w,?^3,0^3)$ satisfies $\ip{a}{b}=\wbar-1$.  So
$R$ has the root $a+b=(2-\wbar;1^6,\wbar,?^3,0^3)$, where the $?$'s are
negated cube roots of~1.

Step 5: {\it $R$ has the roots $(2+\theta;-1^4,1^6,0^3)$ with the
  $0$'s at noncollinear points, the $1$'s on the lines joining them
  and the $-1$'s everywhere else.}  Take $b=(2-\wbar;-1^3,0^3,1^7)$ from
step~4, and try $a=(\w;?,0^2,?,0^2,?^2,0^5)$.  Solving for $?$'s
reveals that $a=(\w;1,0^2,\w,0^2,\w^2,\discretionary{}{}{}0^5)$ satisfies
$\ip{a}{b}=\wbar-1$, so $R$ has the root
$a+b=(2+\theta;0,-1^2,\w,0^2,\discretionary{}{}{}-\wbar^2,1^5)$.

Now we can prove the second claim of the lemma.  If $r$ is in the
first set of roots specified, then $R$ has $r$ by step~4 and
$\p+r=(\thetabar\wbar;1^3,-\w^3,0^7)$ by step~1.  If $r$ is in the
second set of roots, then $R$ has $r$ by step~2 and
$\p+r=(-2+\w;-\w^4,0^3,1^6)$ by step~5.

Finally, we prove that the differences of the $r$'s span $\pperp$.  
We will only need the second batch of roots briefly, so we define $N$
to be the span of the differences of the pairs of roots from the first
batch.  It consists of vectors of the form $(0;\dots)$.  Now we note
that a root from the first batch, minus one from the second, has the
form $(1;\dots)$.  Therefore it suffices to show that $N$ equals the
set $X$ of all vectors $(0;x_1,\dots,x_{13})\in L$ that are orthogonal
to $\p$, which is to say that $x_1+\dots+x_{13}=0$.  We will restrict
attention to the last 13 coordinates.   

Begin by labeling the lines of $P^2\F_3$ by $l_1,\dots,l_{13}$, and
write $r_{ij}$ for the root $(2+\theta;0^3,\wbar^3,-1^7)$ from the
first batch, with the
$0$'s on $l_i$ and the $\wbar$'s on $l_j$.  Then $N$ contains the
vectors $\d_{ij}=-\w(r_{ij}-r_{ji})=(0;1^3,-1^3,0^7)$.  
The span
of the $\d_{ij}$ is easy to understand, because if $i$, $j$, $k$ and $l$
are all distinct, then $\d_{ij}=-\d_{ji}$, $|\d_{ij}|^2=6$,
$\ip{\d_{ij}}{\d_{jk}}=-3$ and $\ip{\d_{ij}}{\d_{kl}}=0$.  It follows
that $N\tensor\C$ admits a coordinate system using 13 coordinates
summing to $0$, in which $\d_{ij}=(\theta,\thetabar,0^{11})$ with
$\theta$ in the $i$th spot and $\thetabar$ in the $j$th.  (One
just checks that the inner products of these vectors, under the
standard pairing, are the same as those of the $\d_{ij}$.)
Write $M$ for the span of the $\d_{ij}$; $L_3(3)$ acts on this
coordinate system by permuting coordinates as it permutes the lines of
$P^2\F_3$.  

Now, $N$ is strictly larger than $M$, because
computation shows that if $i$, $j$ and $k$ are general
lines, then $\ip{r_{ij}-r_{jk}}{\d_{ik}}\notin3\E$.  We can apply
lemma~\ref{lem-unique-enlargement} to both $N$ and $X$ and conclude from the uniqueness
proven there that $N=X$. (We have also shown that $N=X$ admits an
automorphism exchanging the vectors of the form
$(\theta,\thetabar,0^{11})$ with those of the form $(1^3,-1^3,0^7)$.) 
\end{proof}

\begin{proof}[Proof of theorem~\ref{thm-Y555-generates}:]
It remains only to prove $R=\aut L$.  The primitive null vector $\p$
of lemma~\ref{lem-specific-roots-in-R} has Leech type, because
theorem~\ref{thm-Neimeier-lattices} tells us that the complex Leech
lattice is the only Eisenstein Niemeier lattice whose automorphism
group contains $L_3(3)$.  Lemmas~\ref{lem-stabilizer-of-null-vector}
and~\ref{lem-specific-roots-in-R} assure us that $R$ contains the
unipotent radical of the stabilizer of $\p$ ($U$ from the proof of
lemma~\ref{lem-stabilizer-of-null-vector}).  This acts transitively on
the roots $r\in L$ with $\ip{r}{\p}=\theta$, so $R$ contains all their
triflections.  Then the proof of theorem~4.1 of \cite{ch13} shows that
$R$ acts transitively on null vectors of Leech type, so $R$ contains the triflections
in every root having inner product $\theta$ with some null
vector of Leech type.  (These are all the roots of $L$ by
\cite[prop.~4.3]{basak}, but we don't need this.)
The triflections in the point roots obviously
have this property, and those in the line roots do too (by conjugacy).
Therefore $R$ is exactly the group generated by triflections in the
roots with this property, so $R$ is normal in $\aut L$.

Therefore $R$'s intersection with the stabilizer $H$ of $\p$ is normal
in $H$.  Since we already know that $R$ contains $U\normalin H$, $R$
is determined by its image in $H/U\iso 6\suz$.  We also know
(lemma~\ref{lem-Y555-roots-span-and-R-has-L3(3)}) that $R\cap H$ contains $L_3(3)$.  By the simplicity of
$\suz$, $(R\cap H)/U\sset 6\suz$ surjects to $\suz$.  Since $6\suz$
is a perfect central extension of $\suz$, its only subgroup surjecting
to $\suz$ is itself.  Therefore $(R\cap H)/U=6\suz$.  It follows that
$R\cap H$ is all of $H$.  We have shown that $R$ acts transitively on
the primitive null vectors of Leech type, and contains the full
stabilizer of one of them.  So $R=\aut L$.
\end{proof}

\begin{remark}
One  can recover Wilson's $L_3(3)$-invariant
description of the complex Leech lattice $\cll$ (see the end of
\cite{wilson}) by writing down generators for $\pperp$ and then adding
suitable multiples of $\p$ to shift them into $M\tensor_\E\C$, where
$M$ is from the proof of lemma~\ref{lem-specific-roots-in-R}.
\end{remark}

\begin{remark}
We observed that $\rho$ has Leech type.  One can show by patient
calculation that $(\theta;\theta,0^{12})$ has $E_6$ type,
$(\theta;\thetabar,0^{12})$ has $A_2$ type, $(3+\w;1^4,-1^3,0^6)$ has
$D_4$ type, and $(2\theta;\theta^4,0^9)$ has $E_8$ type.  (In the last
case, we specify that the four $\theta$'s are at $4$ points of
$P^2\F_3$ in general position.)
\end{remark}

\end{document}